\pdfoutput=1
\documentclass[12pt]{article}
\usepackage{amssymb,amsmath,amsthm,amscd}
\usepackage{amsfonts}
\usepackage{color}
\usepackage{soul}
\usepackage{hyperref}
\usepackage{microtype}

\makeatletter
\def\@seccntformat#1{\csname the#1\endcsname.\ } 
\makeatother

\def\responseletter{}

\input cyracc.def

\newtheorem{theorem}{Theorem}
\newtheorem{corollary}{Corollary}

\newtheorem{lemma}{Lemma}
\newtheorem{proposition}{Proposition}
\newtheorem{conjecture}{Conjecture}

\theoremstyle{definition}
\newtheorem{definition}{Definition}
\newtheorem{remark}{Remark}

\newcommand\Zech{\mathrm{z}}
\newcommand\MLT{\times} 
\newcommand{\UU}{{\mathcal{L}}}
\newcommand{\FF}{\mathbb{F}}
\newcommand\ZZ{\mathbb{Z}}
\newcommand\BB{\mathcal{B}}
\newcommand\TT{\mathcal{T}}

\newcommand\DD{\mathcal{D}}
\newcommand\VV{\mathcal{V}}
\newcommand\GG{\mathcal{G}}\DeclareMathOperator{\orB}{orb}\DeclareMathOperator{\Aut}{Aut}\DeclareMathOperator{\Cy}{Cyc}
\newcommand\orb[1][G]{\orB_{#1}}
\newcommand\Cos{\mathrm{Cos}}

\newcommand\PG{\mathrm{PG}}

\newcommand\vc[1]{{\bar{#1}}}

\title{\bf \boldmath Triangle decompositions of PG$(n-1,2)$%
\thanks{
This is the accepted version of the manuscript published in Discrete Mathematics 349(1) 2026, 114664(1--13), 
\url{https://doi.org/10.1016/j.disc.2025.114664}. \copyright 2025; available under the CC-BY-NC-ND 4.0 license. 

This research is supported in part by National Natural Science Foundation of China (12471490); the work of Xiaoxiao Li is supported by China Scholarship Council (202306500025); the work
of D.\,S.\,Krotov is supported within the framework
of the state contract of the Sobolev Institute of
Mathematics (Project FWNF-2022-0017).}%
}

\author{
Minjia Shi{$^\dagger$}, Xiaoxiao Li%
\thanks{Key Laboratory of Intelligent Computing Signal
Processing, Ministry of Education, School of Mathematical Sciences, Anhui University, Hefei, 230601, China}%
, Denis S. Krotov%
\thanks{Sobolev Institute of Mathematics, pr. Akademika Koptyuga 4, Novosibirsk, Russia 630090}%
}
\date{}

\newcommand\param[2]{\Delta_{#2}(#1)}
\newcommand\paramg[3]{\Delta^{\mathrm{gd}}_{#3}(#1,#2)}
\newcommand\kgg{{\frac{2^{n}-1}{2^m-1}}}

\newcommand\tri[3]{{\triangle{#1}\raisebox{0.0em}{$#2$}{#3}}}

\begin{document}

\responseletter

\maketitle
\begin{abstract}
We define a triangle design as a partition of the set of lines of a projective space into triangles, where a triangle consists of three pairwise intersecting lines with no common point. A triangle design is balanced if all points are involved in the same number of triangles. We construct balanced triangle designs in PG$(n-1,2)$ for all admissible $n$ (congruent to $1$ modulo $6$) and an infinite class of balanced block-divisible triangle designs. We also prove that the existence of a triangle design in PG$(n-1,2)$ invariant under the action of the Singer cycle group is equivalent to the existence of a partition of $Z_{2^n-1}\backslash\{0\}$ into special $18$-subsets and find such partitions for $n=7$, $13$, $19$.

\end{abstract}

\section{Introduction}

Subspace, or $q$-ary, designs are subspace analogs of classical combinatorial designs,
where the roles of sets and subsets of cardinality~$i$ are played by
spaces and $i$-dimensional subspaces.
For example, balanced incomplete block design and its subspace analog with parameters
$t$-$(v,k,\lambda)$ are defined as a set~$S$ of cardinality~$v$
(subspace of dimension~$v$ over a finite field~$\FF_q$)
and a collection of $k$-subsets of~$S$ (respectively, $k$-dimensional subspaces), called blocks,
such that every $t$-subset (respectively, $t$-dimensional subspace)
is included in exactly $\lambda$ blocks.
Of special interest are designs with $\lambda=1$,
called Steiner systems~\cite{ColMat:Steiner}.
The only known nontrivial parameters of subspace analogs of Steiner systems with $t>1$
are $2$-$(13,3,1)$, $q=2$~\cite{BEOVW:q-Steiner} (see also~\cite{Keevash+:subspace}, where the existence of $t$-$(n,k,\lambda)$ subspace designs is established asymptotically, for sufficiently large~$n$).
Steiner systems of strength $t=2$ can be considered as decompositions
of the complete graph~$K_v$ on~$v$ vertices into complete subgraphs
on~$k$ vertices
such that every edge of~$K_v$ belongs to exactly one subgraph from
the decomposition.
Similar decompositions
of~$K_v$ into isomorphic copies of some other graph~$K$ were
considered in a number of papers, see the survey~\cite{BryElZ:GD}.
A $2$-$(n,K,\lambda)$ design over $\FF_q$,
called a \emph{graph decomposition},
is a collection of graphs (blocks) isomorphic to $K$ with the following properties: the vertex set of every block is a subspace of $\PG(n-1,q)$; every two distinct points of $\PG(n-1,q)$ are adjacent in exactly $\lambda$ blocks, where $\PG(n-1,q)$ is the $(n-1)$-dimensional projective space over $\FF_q$ and $K$ is a simple graph of order $\frac{q^k-1}{q-1}$ for some~$k$, see \cite{BurNakWas:2021}. The definition of a graph decomposition over finite field covers, as a special subcase, the well-known concept of a $2$-$(v,k,\lambda)$ design over finite field, corresponding to the case where $K$ is complete.
In this paper, we consider similar configurations in the subspace case,
decomposing the set of all lines
of a given projective space into copies of one of the simplest substructures,
called a triangle.

We consider the case $q=2$, and for this value
of~$q$, our triangle designs are also
a special case of triangle decompositions
of Steiner triple systems (STS).
Such decompositions were considered in \cite{HorakRosa88}, \cite{Lovegrove:PhD}, \cite{MPZ87};
the results include the existence of
triangle-decomposable STS$(v)$s
for each admissible~$v$~\cite{MPZ87},
the existence of exponential number
of triangle-decomposable STS$(v)$s for infinite number of values of~$v$~\cite{HorakRosa88}, and some recursive constructions that guarantee, in particular,
the existence of triangle decomposition of~$\PG(n-1,2)$ for all $n\equiv 0,1\bmod 6$ \cite[Corollary 7.4.1]{Lovegrove:PhD} (and even for all other~$n$ if we allow one or two lines to be uncovered by triangles).
We aim to construct balanced triangle decompositions,
where ``balanced'' means that all points occur in the same number
of triangles from the decomposition.

Our main results are the following.
We construct balanced triangle designs
in $\PG(n-1,2)$
of all admissible parameters
$n\equiv 1 \bmod 6$
(Theorem~\ref{th:balanced-rec}).
We construct, for every $n\equiv 0\bmod 6$,
a balanced block-divisible triangle design
in~$\PG(n-1,2)$ with
projective $5$-subspaces in the role of blocks
(balanced designs and block-divisible designs
are defined in a standard design-theory manner,
see the next section for formal definitions).
We find symmetric balanced triangle designs
in~$\PG(6,2)$, $\PG(12,2)$, and~$\PG(18,2)$
and develop some theoretical tools,
relating the existence of balanced triangle designs in~$\PG(n-1,2)$
invariant under the action of the Singer cycle group
with the existence of partitions of $\ZZ_{2^{n}-1}\backslash \{0\}$
into $18$-subsets of form $\Gamma(a) \cup \Gamma(b) \cup \Gamma(a+b)$,
where $\Gamma(a)=\{a,\Zech(a),\Zech(-a),-\Zech(-a),-\Zech(a),-a\}$
is the minimal set containing~$a$ and closed under the
operations $x\to -x$ and $x\to\Zech(x)$, where $\Zech(x)$ is Zech's logarithm of~$x$.

\begin{remark}
 The problem of the existence of a partition of
 $\ZZ_{2^{n}-1}\backslash \{0\}$ into $18$-subsets
 of form $\Gamma(a) \cup \Gamma(b) \cup \Gamma(a+b)$
 is very similar to Heffter's difference problem~\cite{Heffter}, solved in \cite{Peltesohn39} and related to the existence of cyclic Steiner triple systems of order~$v$:
 find a partition of
 $\ZZ_{v}\backslash \{0\}$,
 $v\equiv 1 \bmod 6$, or
 $\ZZ_{v}\backslash \{0,\frac{v}{3},\frac{2v}{3}\}$,
 $v\equiv 3 \bmod 6$, into $6$-subsets of form $\Theta(a)\cup\Theta(b)\cup\Theta(a+b)$, where
 $\Theta(a) = \{a,-a\}$.
\end{remark}

The next section contains basic definitions.
The construction of designs with a predefined automorphism group is considered in Section~\ref{sec:comput};
theoretical results can be found in Sections~\ref{sec:inv} and~\ref{sec:6k+1};
the existence of
triangle group divisible designs in~$\PG(11,2)$,
triangle designs in~$\PG(5,2)$,
and balanced triangle designs
in~$\PG(n-1,2)$,
$n=7,13,19$,
are established in
Sections~%
\ref{sec:6_2},
\ref{sec:gd12_6_2},
and~\ref{sec:6k+1},
respectively.
In Section~\ref{s:th},
we describe two theoretical constructions of balanced and unbalanced triangle designs
in~$\PG(n-1,2)$:
one recursive (Section~\ref{s:rec}),
for any admissible~$n$,
and one direct (Section~\ref{sec:6k}), for $n=0\bmod 6$.
Section~\ref{sec:6k} also contains
a construction
of a balanced block-divisible triangle design
in~$\PG(n-1,2)$,
$n\equiv 0\bmod 6$.

\section{Background material}\label{sec:back}
Let $\FF_q$ denote the finite field of order~$q$,
where~$q$ is a prime power.
Let $\FF_q^n$ be the set of all $n$-tuples over~$\FF_q$,
considered as vectors,
with the component-wise addition and multiplication
by a constant.
The definitions below do not specify
a concrete $n$-dimensional vector space~$\VV$,
but since
all $n$-dimensional vector spaces over~$\FF_q$
are isomorphic, one can imply that~$\VV=\FF_q^n$.
However, in constructions, we can vary the representation of~$\VV$, for example,
it can be~$\FF_{q^n}$, or $\FF_{q^m}\times \FF_{q^{n-m}}$,
$0<m<n$, or $\FF_{q^6}^{n/6}$, $n\equiv 0\bmod 6$,
because
$\FF_{q^i}$ is an $i$-dimensional vector
space over~$\FF_{q}$.

To describe the results of computations, we will need a concrete representation
of an extension field of~$\FF_{q}$ (in our paper, $q=2$).
Recall that~$\FF_{q^i}$ can be represented
as the quotient field~$\FF_q[x]/\langle P(x) \rangle$ of the ring of polynomials in~$x$
over~$\FF_q$ modulo a primitive polynomial~$P(x)$ of degree~$i$.
All non-zero elements of~$\FF_{q^i}$ can be written as the~$q^{i}-1$
powers of~$\xi$, where~$\xi $ is a root of~$P(x)$.

Given a vector space $\VV$, the associated projective space $\PG(\VV)$ is
a point--line geometry, where points are one-dimensional subspaces of~$\VV$
and lines correspond to two-dimensional subspaces. If $\VV$ is an $n$-dimensional space over
$\FF_q$, then $\PG(\VV)$ is also denoted $\PG(n-1,q)$, where $n-1$ is called the dimension of the projective space.
In most of this paper, $\VV$ is a vector space over~$\FF_2$, and it is convenient to identify nonzero vectors with points and $3$-sets $\{a,b,a+b\}$ of vectors with lines of~$\PG(\VV)$.
We denote by $\UU(\VV)$ the set of all lines of~$\PG(\VV)$.

\begin{definition}[triangle, triangle design]
Let $\VV$ be an $n$-dimensional space
over~$\FF_{q}$.
The set $\{U_1,U_2,U_3\}$ of three lines
$U_1$, $U_2$, $U_3$ of~$\PG(\VV)$
is called a \emph{triangle}
if $|U_1\cap U_2|=|U_1\cap U_3|=|U_2\cap U_3|=1$ and $|U_1\cap U_2\cap U_3|=0$.
In other words,
a triangle is the set
$\{
\langle a,b \rangle,
\langle b,c \rangle,
\langle c,a \rangle \}$
for some distinct points $a$, $b$, $c$ of~$\PG(\VV)$,
where  {$\langle  x,y \rangle$}
denotes the line through the points $x$ and~$y$.
We will denote such a triangle
by~$\tri{a}{b}{c}$.

A \emph{triangle design} over a finite field or, more precisely, a $\param{n}{q}$-design,
is a set $\TT$ of triangles in~$\PG(\VV)$,
where $\VV$ is an $n$-dimensional vector space over~$\FF_q$,
such that each line of~$\PG(\VV)$ is contained in exactly one triangle in~$\TT$.
\end{definition}

A set~$\GG$ of subspaces of a vector space~$\VV$
such that every non-zero vector in~$\VV$
belongs to exactly one subspace from~$\GG$
is called a \emph{partition} of the vector space~$\VV$.
Partitions of the space into subspaces of the same dimension
are also called \emph{spreads},
and in this paper,
we often use the spread called \emph{Desarguesian},
which consists of the $(q^n-1)/(q^m-1)$
multiplicative cosets of the subfield~$\FF_{q^m}$
of the field~$\FF_{q^n}$, where $m$ divides~$n$.
If~$\FF_{q^n}$ is considered as an $n$-dimensional space over~$\FF_q$,
then each such coset is an $m$-dimensional subspace.
Similarly to the definition of a $q$-analog of group divisible designs (see \cite[Definition 2.1.1]{BKKNW:2019:gdd}), we define the following.
\begin{definition}[group divisible triangle design]\label{TGGD}
Let $m$ be a divisor of~$n$.
A \emph{group divisible triangle design} of order~$n$ over~$\FF_q$,
denoted  $\paramg{n}{m}{q}$,
is a triple $(\VV,\GG,\BB)$, where $\VV$ is an $n$-dimensional vector space over~$\FF_q$, $\GG$ is a partition of~$\VV$ into subspaces (\emph{groups}) of dimension~$m$,
and~$\BB$ is a set of triangles in~$\PG(\VV)$,
that satisfies the following:
\begin{itemize}
 \item triangles from $\BB$ consist of lines that do not lie in any group from~$\GG$,
 \item and any line not lying in any group from~$\GG$ belongs to
 exactly one triangle from~$\BB$.
\end{itemize}
\end{definition}
Simple counting arguments for a $\paramg{n}{m}{q}$-design $(\VV,\GG,\BB)$
yield that there are $ \frac{(q^n-1)(q^n-q^m)}{(q^2-1)(q^2-q)} $
lines that do not lie in one of the subspaces from~$\GG$. Hence,
\begin{equation}\label{eq:BB}
 |\BB|=\frac{1}{3}\cdot\frac{(q^n-1)(q^n-q^m)}{(q^2-1)(q^2-q)}.
\end{equation}
If $q=2$, then the integrality of~\eqref{eq:BB} is equivalent to
holding at least one of the following:
\begin{itemize}
 \item
$m$ (and hence $n$) is even;
 \item
$n$ is divisible by~$6$;
 \item
$n-m$ is divisible by~$6$.
\end{itemize}
It is easy to see that a $\param{n}{q}$-design~$\TT$ is essentially
a $\paramg{n}{1}{q}$-design. So, as the special case $m=1$
of the calculations above, we have
\begin{equation}\label{eq:ness}
 |\TT|=\frac{1}{3}\cdot\frac{(q^n-1)(q^n-q)}{(q^2-1)(q^2-q)},
\end{equation}
and in the binary case (as well as for any~$q$ congruent to~$2$ modulo~$3$), the necessary divisibility condition for the existence of $\param{n}{q}$-designs turns to $$n\equiv 0,1\bmod 6.$$

A set~$\TT$ of triangles is a $\param{n}{q}$-design in~$\PG(\VV)$ if and only if $\TT$ exactly covers~$\UU(\VV)$, i.e.,
$$
\UU(\VV)=\dot{\bigcup_{\{U,V,W\}\in \TT}} \{U,V,W\},
\qquad
\mbox{where $\dot\cup$ denotes the disjoint union}.
$$
We will say that a triangle~$T$ \emph{covers} a point~$\vc x$ if~$\vc x$ belongs to at least one line in~$T$.
\begin{definition}[balanced/unbalanced]\label{balancedesign}
A triangle design or group divisible triangle design $\TT$ in~$\PG(\VV)$ is called \emph{balanced}
if each point of~$\PG(\VV)$ is covered by the same number of triangles from~$\TT$.
If $\TT$ is not balanced, it is \emph{unbalanced}.
\end{definition}
\begin{lemma}\label{kisodd}\label{l:balanced}
Let $\VV$ be an $n$-dimensional space
over~$\FF_{2}$. If $\TT$ is a balanced triangle design in~$\PG(\VV)$,
then $n$ is odd and each point of~$\PG(\VV)$ is covered by $(2^n-2)/3$ triangles from~$\TT$.
\end{lemma}
\begin{proof}
Suppose that each of $2^n-1$
points of~$\PG(\VV)$ is covered by $\mu$ triangles.
Since each triangle covers $6$ different points, we have
$(2^n-1)\mu=6|\TT|$. Therefore, $|\TT|={(2^n-1)\mu}/{6}$.
On the other hand, $|\TT|={(2^n-1)(2^n-2)}/{18}$, see~\eqref{eq:ness}. By equating,
we find that $\mu$ equals $(2^n-2)/3$, which is integer if and only if $n$ is odd.
\end{proof}

\section{Triangle designs with a prescribed automorphism group}\label{sec:comput}

In this section, we consider the construction
of triangle designs with some prescribed automorphism
groups. While the concrete existence results of this section
(Theorems~\ref{th:6_2}, \ref{th:gd12-6}, \ref{th:n7n13n19})
are computational, Sections~\ref{sec:inv} and~\ref{sec:6k+1} include also important
theoretical conclusions that reduce the existence of triangle designs invariant under certain automorphism groups
to some numerical configurations, which are conjectured
to exist for all parameters satisfying the straightforward divisibility conditions.

Finding a triangle design is obviously a special case
of the \emph{exact covering} problem:
we need to find a partition of the given set
(in our case, the set of all lines)
into subsets from the given family
(in our case, the family of all triangles).
For small parameters, this problem can be solved
with existing software, e.g., \cite{KasPot08}.
However, the size of the exact covering problem
grows rapidly with increasing the parameters of the original task,
and a standard approach is to reduce the task
with prescribing some symmetry to a covering to find.
If we imply that the covering is invariant under the action
of some prescribed automorphism group,
then the original task can be reduced to the problem
of covering the orbits of elements by the orbits of subsets.
This approach is well described in the famous paper
of Kramer and Mesner~\cite{KraMes:76},
and there is no need to explain it here.
We only note that the automorphisms we prescribe,
the Singer cycle and the Frobenius automorphism,
are the most popular in searching different kinds of
designs over subspaces of a finite vector space,
see, e.g., \cite{BEOVW:q-Steiner}, \cite{BKW:2018subspace}, \cite{BKOW:2014}.

Let $G$ be an automorphism group of the $n$-dimensional
vector space~$\VV$ over~$\FF_q$.
The $G$-orbit of a line~$L$ is denoted by
$\orb (L) :=\{\pi(L):\pi\in G\}\subseteq \UU(\VV)$.
Further, we find convenient to introduce the following
concept.

\begin{definition}[triangle-orbit]
If $\{U,V,W\}$ is a triangle and
 $\orb(U)\ne\orb(V)\ne\orb(W)\ne\orb(U)$,
 then we call the set $\{\orb(U),\orb(V),\orb(W)\}$
 a \emph{triangle-orbit}.
\end{definition}

In the subsections below, we consider several computational results,
establishing the existence of
$\param{6}{2}$, balanced $\paramg{12}{6}{2}$,
and balanced $\param{n}{2}$, $n=7,13,19$,
and some theoretical results concerning triangle designs invariant under
the action of $\FF_{2^n}^\MLT$ (Section~\ref{sec:inv})
and, additionally, $\Aut(\FF_{2^n})$  (Section~\ref{sec:6k+1}).

\subsection
[The existence of an unbalanced (6)2-triangle-designs]
{The existence of an unbalanced \boldmath $\param{6}{2}$-designs}
\label{sec:6_2}
The existence of a  $\param{6}{2}$-design
(by Lemma~\ref{l:balanced}, it is unbalanced)
is not only
a stand-alone result,
but also used as a brick
in theoretical constructions, both
recursive (Section~\ref{s:rec}) and
direct (Section~\ref{sec:6k}).
\begin{theorem}\label{th:6_2}
 There is a $\param{6}{2}$-design with an automorphism group isomorphic to~$\FF_{32}^\MLT$.
\end{theorem}

\begin{proof}
Taking into account the isomorphism between vector spaces,
we consider $\VV = \FF_2\times\FF_{2^5}$,
where
$\FF_{2^5} =
\FF_2[x]
/
\langle x^5+ x^2 +1 \rangle
$,
as a $6$-dimensional vector
space over~$\FF_2$.
For $b\in \FF^\MLT_{2^5}=\FF_{2^5}\setminus \{0\}$,
the mapping~$ \mu_b: \VV \to \VV $ defined as
$$
\mu_b((a,y)) = (a,b\cdot y),
\quad a\in \FF_2,\ y\in \FF_{2^5},
$$
is an automorphism of the vector space~$\VV$.
Let $G=\{\mu_b:\ b\in \FF^\MLT_{2^5}\}$. Then $G$ is a group of automorphisms of~$\VV$ of order~$2^5-1=31$.
There are~$651$ lines in~$\PG(\VV)$,
i.e., $|\UU(\VV)|=651=21\cdot 31$. We observe that $G$ does not fix any line, i.e., acts semiregularly on~$\UU(\VV)$.
Hence, $G$ partitions~$\UU(\VV)$ into $|\UU(\VV)|/|G|=21$ orbits.
%
Solving the corresponding exact covering problem,
we find a $\param{6}{2}$-design consisting of $7$ orbits of triangles
with the following representatives:
$$
  [0^{3},0^{9},1^{27}]
,\ [0^{4},0^{5},1^{0}]
,\ [0^{4},0^{13},1^{15}]
,\ [0^{27},0^{29},1^{10}]
,\ [0^{3},1^{3},1^{7}]
,\ [0^{5},0^{28},1^{4}]
,\ [1^{13},1^{18},1^{25}],
$$
where
$[a^\alpha,b^\beta,c^\gamma]$ encodes the triangle
$ \tri{(a, \xi^\alpha)}{(b, \xi^\beta )} {(c, \xi^\gamma)}
$.
\end{proof}

\subsection[Triangle designs invariant under the Singer cycle group]{Triangle designs invariant under the Singer cycle group}\label{sec:inv}

In this section, we consider group divisible triangle designs
(and triangle designs, as a special case) in $\FF_2^{n}$
that are invariant under the action of the \emph{Singer cycle},
the cyclic group of order $2^n-1$ corresponding to the multiplicative group
$ \FF_{2^{n}}^\MLT $ of the field~$ \FF_{2^{n}}$.

\begin{lemma}\label{l:orbB}
 Assume that $(\VV,\GG,\BB)$ is a $\paramg{n}{m}{2}$-design,
 where $\VV = \FF_{2^n}$,
 $\GG$ is the Desarguesian spread
 and $\BB$ is invariant under the action of the Singer cycle group
 $G =  \FF_{2^{n}}^\MLT $.
 If $\{B_1,B_2,B_3\}$ is a triangle from $\BB$, then the orbits
        $\orb(B_1)$, $\orb(B_2)$, $\orb(B_3)$ are distinct and their cardinalities
        are equal to~$2^n-1$.
\end{lemma}
\begin{proof}
Trivially, the claim is equivalent to the following:
\begin{itemize}
 \item[(*)] \em For a triangle
$\tri{a}{b}{c}$
 from~$\BB$,
 all $3(2^{n}-1)$ lines $\langle \gamma a,\gamma b \rangle$,  $\langle \gamma b,\gamma c \rangle$,  $\langle \gamma c,\gamma a \rangle$, $\gamma\in\FF_{2^{n}}^\MLT$, are distinct.
\end{itemize}
To prove (*), seeking a contradiction,
suppose that $\langle \gamma a,\gamma b \rangle$ coincides with
$\langle \gamma' a,\gamma' b \rangle$, $\langle \gamma' b,\gamma' c \rangle$,
or $\langle \gamma' c,\gamma' a \rangle$, where
$\gamma\ne\gamma'$.
Since the triangles
$$\tri {(\gamma a)}{(\gamma b)}{(\gamma c)}\quad\mbox{and}\quad\tri {(\gamma' a)}{(\gamma' b)}{(\gamma' c)}$$
are both in~$\BB$ and cover the same line,
$\langle \gamma a,\gamma b \rangle$,
they are the same.
It follows that $\{\gamma' a,\gamma' b,\gamma' c\} = \{\gamma a,\gamma  b,\gamma c\}$
and hence $\gamma' (a+b+c) = \gamma(a+b+c)$,
which implies $a+b+c=0$. Then we have
$\langle a,b \rangle = \langle b,c \rangle = \langle c,a \rangle$,
which contradicts the definition of a triangle.
\end{proof}

\begin{corollary}\label{c:mod6}
Under the hypothesis of Lemma~\ref{l:orbB}, it holds $n-m \equiv 0 \bmod 6$.
\end{corollary}
\begin{proof}
From \eqref{eq:BB}, $q=2$,
we have $|\BB|=(2^n-1)(2^n-2^m)/18$.
By Lemma~\ref{l:orbB}, $|\BB|$ is divisible by~$2^n-1$. It follows that $(2^n-2^m)/18$
is an integer, in particular, $2^{n-m}-1$ is divisible by~$9$.
This happens if and only if $n-m \equiv 0 \bmod 6$.
\end{proof}

In the rest of this section, $\VV = \FF_{2^{n}}= \FF_2[x]/\langle P(x) \rangle$,
where $n=sm$, $m\ge 1$, $s\ge 2$, $n-m \equiv 0 \bmod 6$.
Let $\xi$ be a primitive element of~$\FF_{2^{n}}$ and a root of~$P(x)$.
Then $\xi^\kgg$ is a primitive element of~$\FF_{2^{m}}$,
a subfield of~$\FF_{2^{n}}$.
The multiplicative group~$\FF^\MLT_{2^m}$ is a subgroup of~$\FF^\MLT_{2^{n}}$.
The number of cosets of~$\FF^\MLT_{2^m}$ in~$\FF^\MLT_{2^{n}}$
is~$\kgg$.
Actually, $\FF^\MLT_{2^m}$ is decomposed into cosets of~$\FF^\MLT_{2^{n}}$ as follows:
$$\FF^\MLT_{2^{n}}=\FF^\MLT_{2^m}\cup \xi \FF^\MLT_{2^m}\cup\cdots\cup \xi^{\kgg-1}\FF^\MLT_{2^m}.$$
If $U=\langle a,b\rangle$ is a line of~$\PG(\VV)$, then
either
\begin{enumerate}
\item [(1)] $U$ is a subset of $\xi^i\FF_{2^m}=\xi^i\FF^\MLT_{2^m}\cup \{0\}$
   for some $i \in \{0,1,\ldots,\kgg-1 \}$, or
\item [(2)]
$a\in \xi^i\FF_{2^m}$,
$b\in \xi^j\FF_{2^m}$,
$a+b\in \xi^k\FF_{2^m}$
for distinct $i,j,k \in \{0,1,\ldots,\kgg-1 \}$.
\end{enumerate}
Let $\GG=\{\xi^i\FF_{2^m}:\ 0\leq i\leq \kgg-1 \}$.
Then $\GG$ is a partition of~$\VV$
into subspaces of dimension~$m$.
Let
\begin{align*}
   \UU_i&=\big\{L\in \UU(\VV):\ L\subseteq \xi^i\FF_{2^m}\big\},\qquad \textstyle i=0,\ldots,\kgg-1;\\
   \overline{\UU}&=\UU(\VV)\setminus \big(\UU_0 \cup \ldots \cup \UU_{\kgg-1}\big).
\end{align*}
To find a symmetric $\paramg{n}{m}{2}$-design,
we need to partition $\overline{\UU}$
into triangle-orbits.

If $a\in \FF^\MLT_{2^{n}}$, then $\pi_a$ defined as
$\pi_a(y)=ay$ is an automorphism of the vector space~$\VV$.
Let $G=\{\pi_a:\ a\in \FF^\MLT_{2^{n}}\}$.
Then $G$ is a group of automorphisms of~$\VV$
of order~$2^{n}-1$.
Moreover,
$G$ acts semiregularly on~$\overline{\UU}$
(but not on $\UU(\VV)$ if $m$ is even, because $\pi_{\xi^{(2^{n}-1)/3}}$
fixes $\{0,1,\xi^{(2^{n}-1)/3},\xi^{2(2^{n}-1)/3}\} \subset \UU_0$).

Let $K=\{1,2,\ldots,2^{n}-2\}$.
Since $\xi^{2^{n}-1}=\xi^0=1$,
we consider~$K$ as a subset of~
$\ZZ_{2^{n}-1} =
\ZZ / (2^{n}-1) \ZZ$, with operations modulo~$2^{n}-1$.
For $k\in K$, we denote
by~$\Zech(k)$ the \emph{Zech logarithm}
$\log_\xi(1+\xi^k)$, i.e., $1+\xi^k=\xi^{\Zech(k)}$.
Let $\overline{K}=K\setminus \{k\in K: \kgg|k\}$. For each $k\in \overline{K}$, denote
\begin{equation}
 \label{eq:Gamma}
\Gamma(k)=\{k,\Zech(k),\Zech(-k),-\Zech(-k),-\Zech(k),-k\}.
\end{equation}
\begin{proposition}\label{prop1}
The following assertions are true:
\begin{itemize}
 \item [\rm(i)] the size of~$\Gamma(k)$ is~$6$ for each~$k$ in~$K$ such that $3k\ne 0$
 (in particular, for each~$k$ in~$\overline{K}$);
 \item [\rm(ii)] for each~$k$ in~$K$ and~$k'$ in~$\Gamma(k)$,
 it holds $\Gamma(k')=\Gamma(k)$;
 \item [\rm(iii)] if $k\in \overline{K}$, then $\Gamma(k) \subset \overline{K}$.
\end{itemize}
\end{proposition}
\begin{proof}
Since,
\begin{equation}
\xi^{\Zech(-\Zech(k))}
= 1+\frac{1}{1+\xi^{k}}
= \frac{1}{\xi^{-k}+1} =
\xi^{-\Zech(-k)}
= \frac{\xi^k}{1+\xi^k}
=\xi^{k-\Zech(k)},
\end{equation}
we have
$\Zech(-\Zech(k))=-\Zech(-k) = k-\Zech(k)$.
Additionally, we have $-(-k)=k$ and $\Zech(\Zech(k))=k$.
With all these identities,
it is easy to see that by applying~$-$ or~$\Zech()$
to any element of~$\Gamma(k)$,
we get again an element of~$\Gamma(k)$.
So, (ii) is true.

By (ii), every two elements~$a$ and~$b$ of~$\Gamma(k)$ are in one of the following relations:
$a=b$, $a=-b$, $a=\Zech(b)$, $a=-\Zech(b)$, $a=\Zech(-b)$, $a=-\Zech(-b)$.
If $a=-b$ or $a=\Zech(b)$, then obviously $a\ne b$.
If $a=-\Zech(-b)$, then
$\xi ^a = \frac{1}{1+1/\xi ^b} =\frac{\xi ^b}{1+\xi ^b} \ne \xi ^b$ and hence $a\ne b$.
To prove~(i), it remains to exclude the possibility $b=-\Zech(b)$ (the second case $b=\Zech(-b)$ is similar), which leads to the equation
\begin{equation}\label{eq:2bb}
\xi^{2b} + \xi^b +1 =0.
\end{equation}
Multiplying by $\xi^b-1$,
we get $\xi^{3b}=1$.
In this case, $3b=0$ and
$\Gamma(k)=\Gamma(b)=\{b,-b\}$.
This completes the proof of~(i).

To see (iii), we observe that if $k'$
from~$\Gamma(k)$ is not from $\overline K$,
then $\xi^{k'}$ belongs to the subfield~$\FF_{2^m}$.
In this case, since $k \in \Gamma(k')$ by~(ii),
there also hold $\xi^{k}\in \FF_{2^m}$
and~$k\not\in \overline{K}$.
\end{proof}

\begin{proposition}\label{prop2}
For each $U\in {\UU(\VV)}$, there is $k\in {K}$ such that
\begin{equation}
 \label{eq:repr}
\orb(U)
= \orb\langle\xi^k,\xi^{\Zech(k)}\rangle
= \orb\langle \xi^{-k},\xi^{\Zech(-k)}\rangle
= \orb\langle \xi^{-\Zech(-k)},\xi^{-\Zech(k)}\rangle.
\end{equation}
Moreover, $\langle\xi^k,\xi^{\Zech(k)}\rangle$,
$\langle \xi^{-k},\xi^{\Zech(-k)}\rangle$, and
$\langle \xi^{-\Zech(-k)},\xi^{-\Zech(k)}\rangle$
are the only representatives of~$\orb(U)$
that contain~$1$.
If $U \in \overline\UU$, then $k\in \overline K$.
\end{proposition}
\begin{proof}
Suppose $L=\{a,b,a+b\}$ is a line from $\overline\UU$. It is easy to see that $b^{-1}L=\{1,\frac{a}{b},\frac{a+b}{b}\}$ is a line in~$\orb(L)$.
There is~$i\in\{0,\ldots,\frac{2^n-1}{2^m-1}-1\}$ such that $\frac{a}{b}\in \xi^iF^\MLT_{2^{m}}$.
So, $\frac{a}{b}=\xi^k=\xi^{i+ \kgg j}$ for some~$j$. Since $1+\xi^k=\xi^{\Zech(k)}$, we see that
$\{1,\xi^k,\xi^{\Zech(k)}\}
=
\langle 1
, \xi^k
\rangle$
is a line in $\orb(L)$. Similarly,
$
\{ 1, \xi^{-k},\xi^{\Zech(-k)}\}
=
\langle 1, \xi^{-k}\rangle
$
and
$
\{ 1,\xi^{-\Zech(k)},\xi^{\Zech(-\Zech(k))}\}
=
\langle 1,\xi^{-\Zech(k)}\rangle
$
belong to~$\orb(L)$.
To complete the proof of~\eqref{eq:repr},
it remains to show that $\Zech(-\Zech(k)) = -\Zech(-k)$:
$$
\xi^{\Zech(-\Zech(k))}
= 1+\frac{1}{1+\xi^{k}}
= \frac{1}{\xi^{-k}+1} =
\xi^{-\Zech(-k)}.
$$

Next, for any
$L=\{a,b,a+b\}\in \UU(\VV)$,
the $G$-orbit of~$L$
has at most~$3$
representatives
containing~$1$:
$\{1,\frac ba, 1+\frac ba\}$,
$\{\frac ab,1,\frac ab+1\}$,
and
$\{\frac{a}{a+b},\frac{b}{a+b}, 1\}$.
All of them are listed in~\eqref{eq:repr}, which proves
the second claim.

Finally,
if $k$ is not in~$\overline K$,
then it is divisible by~$\kgg$,
$\langle 1, \xi^{-k}\rangle$ is in~$\UU_0$,
every element of the orbit
is in~$\UU_i$ for some $i$,
and $L$ is not in~$\overline\UU$.
\end{proof}

\begin{proposition}\label{prop3}
Let $\Gamma(k_i)$ represent $\orb(L_i)$, $i=1,2,3$, in the sense of~\eqref{eq:repr}. {Suppose all three orbits are distinct.}
The triple $\{\orb(L_1),\orb(L_2),\orb(L_3)\}$ is a triangle-orbit if and only if
$s_1 + s_2 + s_3 =0$
for some
$s_1\in \Gamma(k_1)$,
$s_2\in \Gamma(k_2)$,
$s_3\in \Gamma(k_3)$.
\end{proposition}
\begin{proof}
Suppose that $\{\orb(L_1),\orb(L_2),\orb(L_3)\}$ is a triangle-orbit. Then  there are $U\in \orb(L_1)$, $V\in \orb(L_2)$, and $W\in \orb(L_3)$ such that $\{U,V,W\}$ is a triangle. Therefore, $U=\{a,b,a+b\}$, $V=\{a,c,a+c\}$ and $W=\{b,c,b+c\}$ for some $a$, $b$,~$c$.
Let $\frac{a}{b}=\xi^{s_1}$, $\frac{c}{a}=\xi^{s_2}$,
and $\frac{b}{c}=\xi^{s_3}$. In particular,
we see $s_1+s_2+s_3=0$.
By Proposition~\ref{prop2}, we have
$\Gamma(k_i)=\Gamma(s_i) $, $i=1,2,3$.

Conversely, assume that
$\orb(L_i)$ contains $\langle \xi^{k_i},\xi^{\Zech(k_i)} \rangle$ for each $i$ in $\{1,2,3\}$ and that $\xi^{s_1}\xi^{s_2}\xi^{s_3}=1$
for some
$s_1\in \Gamma(k_1)$,
$s_2\in \Gamma(k_2)$,
$s_3\in \Gamma(k_3)$.
By Propositions~\ref{prop1} and~\ref{prop2},
$\{ 1,\xi^{s_i},\xi^{\Zech(s_i)}\}$ is a line from~$\orb(L_i)$, $i=1,2,3$.
Denote
$a=\xi^{s_1}$,
$c=\xi^{s_1}\xi^{s_2}=\xi^{-s_3}$,
$U=\{ 1,\xi^{s_1},\xi^{\Zech(s_1)}\}=\{ 1,a,a+1\}$, $V=\xi^{s_1}\{ 1,\xi^{s_2},\xi^{\Zech(s_2)}\}=\{ a,c,a+c\}$, and
$W= \{ 1,\xi^{-s_3},\xi^{\Zech(-s_3)}\}=\{ 1,c, c+1\}$.
If $a+1=c$, then $U=W$ and hence
$\orb(L_1)=\orb(L_3)$, a contradiction.
Therefore, $\{U,V,W\}$ is a triangle and
$\{\orb(L_1),\orb(L_2),\orb(L_3)\}$
is a triangle-orbit.
\end{proof}

\begin{theorem}\label{th:mult-suff}
If $m$ is a divisor of $n$,
then the following assertions are equivalent:
\begin{itemize}
\item[\rm(i)] there exists
a (balanced) $\paramg{n}{m}{2}$-design $(\VV,\GG,\BB)$,
where $\VV = \FF_{2^{n}}$,
$\GG$ is the Desarguesian spread (the set of multiplicative cosets
of the subfield $\FF_{2^{m}}$ of $\FF_{2^{n}}$),
and $\BB$ is invariant under the action of the multiplicative group
$G=\FF^\MLT_{2^{n}} ;$
 \item[\rm(ii)] there exists a partition of
$ \ZZ_{2^{n}-1} \setminus \kgg\ZZ_{2^{n}-1} $
into $18$-subsets of form
$$
\Gamma(k_1)\,\dot\cup\,\Gamma(k_2)\,\dot\cup\,\Gamma(k_3),
\quad
 k_1+k_2+k_3=0,
$$
or, equivalently, a partition of
$ \{\Gamma(k):\ k\in \ZZ_{2^{n}-1} \setminus \kgg\ZZ_{2^{n}-1} \} $
into $3$-subsets of form
$$
\{\Gamma(k_1),\Gamma(k_2),\Gamma(k_3)\},
\quad
k_1+k_2+k_3=0.
$$
\end{itemize}
\end{theorem}

Note that $\Gamma(\cdot)$ and hence the concrete partitions in (ii) depend on the concrete
choice of the primitive element~$\xi$, but the existence of such a partition does not depend
on this choice.


\begin{proof}[Proof of Theorem~\ref{th:mult-suff}]
We first note that each of (i), (ii) implies $n-m \equiv 0 \bmod 6$,
so we can assume it holds.

By Lemma~\ref{l:orbB},
for each triangle $\tri{a}{b}{c}$ in~$\BB$, the triple
$\{ \orB_G\langle  a, b \rangle,$\ $  \orB_G\langle  b, c \rangle,$\ $\orB_G\langle  c, a \rangle \}$
is a triangle-orbit. Hence,
a $\paramg{n}{m}{2}$-design $(\VV,\GG,\BB)$ invariant under the action of~$G$ is equivalent to
a partition of the set of $G$-orbits in~$\overline\UU$ into triangle-orbits.

By Proposition~\ref{prop1},
$\overline{K}$ is the disjoint union of $6$-subsets of form~$\Gamma(k)$.
By Proposition~\ref{prop2},
\eqref{eq:repr} gives a one-to-one correspondence between those $6$-subsets
and $G$-orbits in $\overline\UU$.
By Proposition~\ref{prop3}, three such $G$-orbits form a triangle-orbit
if and only if the corresponding $6$-subsets of $\overline{K}$ are
of form
$\Gamma(k_1)$, $\Gamma(k_2)$, $\Gamma(k_3)$ for some $k_1$, $k_2$, $k_3$
such that $k_1+k_2+k_3=0$.

It follows that (i) and (ii) are equivalent.
\end{proof}

\subsection
[The existence of balanced (12,6,2)- and (6,2,2)-group-divisible triangle design]
{The existence of balanced \boldmath $\paramg{12}{6}{2}$- and $\paramg{6}{2}{2}$-designs}
\label{sec:gd12_6_2}
In this section, $\VV = \FF_{2^{12}}= \FF_2[x] / \langle x^{12}+x^7+x^6+x^5+x^3+x+1\rangle
$.
By Proposition~\ref{prop3},
the covering of~$\overline\UU$ with $224$ disjoint triangle-orbits is equivalent
to a partition of the set $\{\Gamma(a):\ a\in\overline K\}$
into $3$-subsets $\{\Gamma(a),\Gamma(b),\Gamma(c)\}$ such that $a+b+c=0$.
Such partition is constructed by solving the corresponding exact covering problem.
Representatives of triangles
of the form
$\tri{\xi^0}{\xi^i}{\xi^j}$
are listed in
{Table~\ref{t:12}} as the pairs $(i,j)$.
Finally, we get the following.
\begin{theorem}\label{th:gd12-6}
There exists a balanced $\paramg{12}{6}{2}$-design $(\VV,\GG,\BB)$,
where $ \VV=\FF_{2^{12}}$,   $\GG$ is the set of all $65$
$1$-dimensional subspaces of~$\VV$ over~$\FF_{2^6}$
(i.e., the Desarguesian spread; each group of~$\GG$ is a $6$-dimensional vector space over $\FF_2$),
and $\BB$ is invariant under the action of~$\FF_{2^{12}}^\MLT$.
\end{theorem}

The design from the following theorem is used in the construction in Section~\ref{s:rec}.
\begin{theorem}\label{th:gd6-2}
There exists a balanced $\paramg{6}{2}{2}$-design $(\VV,\GG,\BB)$,
where $ \VV=\FF_{2^{6}}$, $\GG$ is the set of all $21$
$1$-dimensional subspaces of~$\VV$ over~$\FF_{2^2}$
(i.e., the Desarguesian spread; each group of~$\GG$ is a $2$-dimensional vector space over $\FF_2$).
\end{theorem}
\begin{proof}
    An example of the triangle set of such a design
is
\begin{multline*}
   \BB= \big\{ \,
    \tri{\xi^{i+3l}}{\xi^{j+3l}} {\xi^{k+3l}}
: \\(i,j,k)\in \{(0, 5, 61),\ (0, 8, 48),\ (0, 9, 58),\ (0, 13, 56),\ (0, 14, 46),\ (0, 17, 20),\\ (0, 28, 31),\ (0, 43, 52),\ (1, 16, 41),\ (2, 11, 26)\},\ l\in\{0,\ldots,20\}\,\big\},
\end{multline*}
where $\xi$ is a primitive root of
$P(z)=z^6 + z^4 + z^3 + z + 1$
in $\FF_{2^6} = \FF_2[z] / \langle P(z) \rangle $.
\end{proof}

\begin{table}
\small
(3,1861),
 (7,1775),
 (15,1342),
 (20,1605),
 (31,321),
 (33,186),
 (34,2042),
 (35,537),
 (36,1492),
 (38,258),
 (41,797),
 (42,1747),
 (43,703),
 (44,171),
 (49,144),
 (51,863),
 (56,556),
 (57,1778),
 (58,1937),
 (68,1054),
 (69,724),
 (75,1394),
 (78,1707),
 (80,344),
 (84,1183),
 (85,250),
 (89,1847),
 (94,775),
 (98,1284),
 (101,834),
 (103,331),
 (113,1081),
 (114,811),
 (117,1158),
 (119,928),
 (120,1759),
 (123,1982),
 (124,1409),
 (134,1025),
 (136,1878),
 (137,1703),
 (141,1203),
 (142,669),
 (143,946),
 (146,2001),
 (147,1571),
 (148,433),
 (151,1044),
 (155,437),
 (157,1073),
 (161,1039),
 (162,1493),
 (167,1839),
 (168,1258),
 (173,1123),
 (174,1415),
 (176,1888),
 (177,1688),
 (179,1706),
 (184,1544),
 (188,518),
 (189,506),
 (190,1458),
 (197,1793),
 (209,982),
 (210,1814),
 (211,2013),
 (212,1199),
 (216,1316),
 (217,1883),
 (225,1289),
 (226,1449),
 (230,1230),
 (232,848),
 (233,1457),
 (234,1632),
 (237,619),
 (240,1951),
 (243,787),
 (247,1986),
 (249,1654),
 (255,727),
 (259,1938),
 (261,1057),
 (265,1298),
 (266,712),
 (267,1517),
 (268,842),
 (273,1212),
 (280,1220),
 (281,1995),
 (287,927),
 (293,1460),
 (295,680),
 (301,608),
 (303,1161),
 (306,1427),
 (312,631),
 (316,1752),
 (335,1860),
 (336,1113),
 (337,1640),
 (340,2136),
 (347,1087),
 (351,793),
 (357,944),
 (360,839),
 (365,1601),
 (366,1811),
 (367,1618),
 (368,1524),
 (370,1249),
 (371,1055),
 (379,900),
 (380,869),
 (384,1397),
 (386,1350),
 (389,1083),
 (391,2085),
 (395,976),
 (403,1326),
 (405,2014),
 (407,1549),
 (410,1574),
 (416,1801),
 (418,1763),
 (419,1884),
 (425,1254),
 (428,1969),
 (448,1239),
 (454,1638),
 (457,1622),
 (461,2195),
 (478,1502),
 (480,1256),
 (486,1012),
 (487,1480),
 (488,1579),
 (493,1538),
 (494,1454),
 (495,1202),
 (496,1367),
 (501,2052),
 (503,1209),
 (511,1393),
 (517,1736),
 (519,1163),
 (523,2153),
 (529,1435),
 (541,2038),
 (552,1297),
 (559,2018),
 (566,1737),
 (570,2033),
 (586,1378),
 (588,1748),
 (591,1540),
 (606,1894),
 (611,1425),
 (626,2260),
 (630,1392),
 (632,2297),
 (639,2173),
 (643,2090),
 (645,2154),
 (646,1819),
 (647,2186),
 (651,1580),
 (654,1491),
 (659,2046),
 (661,1789),
 (664,1477),
 (671,2307),
 (675,2243),
 (679,1652),
 (682,1717),
 (686,1504),
 (690,1983),
 (692,1816),
 (698,1815),
 (709,1453),
 (711,2272),
 (716,1542),
 (728,2111),
 (730,2205),
 (736,2393),
 (737,1764),
 (743,2004),
 (754,2164),
 (766,1954),
 (788,2065),
 (790,2254),
 (794,2026),
 (795,2070),
 (802,1743),
 (804,1949),
 (807,1662),
 (824,1709),
 (833,2376),
 (844,2166),
 (846,2435),
 (854,2189),
 (860,2143),
 (883,1962),
 (902,2161),
 (907,2391),
 (926,2411),
 (937,2183),
 (942,2415),
 (983,2313),
 (990,2479),
 (991,2344),
 (1004,2181),
 (1016,2349),
 (1048,2425),
 (1049,2283),
 (1061,2530),
 (1139,2382),
 (1191,2498),
 (1194,2616),
 (1197,2573),
 (1231,2525),
 (1257,2627),
 (1278,2574)
\caption{The pairs $(i,j)$ encoding the representatives $
\triangle{1}{\xi^i}{\xi^j}
$ of triangles of $\paramg{12}{6}{2}$.}
\label{t:12}
\end{table}

\subsection
[Balanced   triangle designs in PG(6k,2)]
{Balanced   triangle designs in \boldmath PG$(6k,2)$.}
\label{sec:6k+1}

The special case $m=1$ of Theorem~\ref{th:mult-suff} is the following.

\begin{corollary}\label{c:mult-suff}
The following assertions are equivalent:
\begin{itemize}
\item[\rm(i')] in $\PG(\FF_{2^{n}})$, there exists
a (balanced) $\param{n}{2}$-design invariant under the action of the multiplicative group
$\FF^\MLT_{2^{n}} ;$
 \item[\rm(ii')] there exists a partition of
$ \ZZ_{2^{n}-1} \setminus \{0\} $
into $18$-subsets of form
$$
\Gamma(a)\,\dot\cup\,\Gamma(b)\,\dot\cup\,\Gamma(c)
\quad
\mbox{where } a+b+c=0.
$$
\end{itemize}
\end{corollary}
Clearly, for (ii') to hold, $2^n-2$ should be divisible by~$18$,
i.e., $n\equiv 1 \bmod 6$.
\begin{conjecture}\label{conj1}
Assertion {\rm (ii')} holds for each $n=6k+1$.
\end{conjecture}

Conjecture~\ref{conj1} with Corollary~\ref{c:mult-suff} implies
that for each $n=6k+1$ there exists
a $\param{n}{2}$-design
invariant under the action of~$\FF_{2^n}^\MLT$.
Clearly, (ii') can be treated
as an instance of the exact covering problem.
For small~$n$ ($n=7$),
it can be solved computationally.
To solve larger values of~$n$,
we can try to require more symmetries.
Below, we consider triangle designs
that are invariant not only under
the multiplicative group,
but also under the group~$\Aut(\FF_{2^n})$
of automorphisms of the field~$\FF_{2^n}$.
It is known that
$\Aut(\FF_{2^n}) = \{ \phi^0,\phi^1,\ldots,\phi^{n-1}\}$, where $\phi(y)=y^2$,
the Frobenius automorphism (see e.g.~\cite[Ch.\,1, \S4]{Weil:95}).
The partition of
$\FF_{2^n}^\MLT = \{\xi^0,\ldots,\xi^{2^n-2}\}$
into orbits under~$\Aut(\FF_{2^n})$ corresponds
to the partition of
$\ZZ_{2^n-1} = \{0,\ldots,{2^n-2}\}$
into orbits under the multiplication by~$2$,
called \emph{$2$-cyclotomic classes}.
For $k \in \ZZ_{2^n-1}$,
denote by $\Cy(k)$ the $2$-cyclotomic class
containing~$k$.

\begin{proposition}\label{p:cosets2}
Assume that
$k,k'\in\ZZ_{2^n-1} \setminus\{0\}$.
\begin{itemize}
 \item[\rm(i)]  If $\Cy(k)=\Cy(k')$, then
 $\Cy(-k)=\Cy(-k')$ and $\Cy(\Zech(k))=\Cy(\Zech(k'))$.
 \item[\rm(ii)]
  If $n$ is odd, then
 $\Cy(k)\ne\Cy(-k)$,
 $\Cy(k)\ne\Cy(\Zech(k))$,
 and
 $\Cy(k)\ne\Cy(-\Zech(-k))$.
 \item[\rm(iii)] If $n$ is not divisible by~$3$
 and $3k\ne 0$, then $\Cy(-\Zech(k)
)\ne\Cy(k)\ne\Cy(\Zech(-k))$.
\end{itemize}
\end{proposition}
\begin{proof}
(i)
Assume $\Cy(k)=\Cy(k')$.
In this case,
there is $i\in \{0,\ldots,{2^n-2}\}$ such that $k'=2^ik$.
Hence, $-k'=2^i(-k)$ and $\Cy(-k)=\Cy(-k')$.

Since $\xi^{\Zech(k')}=1+\xi^{k'}=1+\xi^{2^ik}=(1+\xi^k)^{2^i}=\xi^{2^i \Zech(k)}$, we have $\Zech(k')\equiv2^i \Zech(k) \bmod 2^{n}-1$. Therefore, $\Cy(\Zech(k))=\Cy(\Zech(k'))$.

(ii)
Let $\sigma(k)$ denote $-k$, $\Zech(k)$, or $-\Zech(-k)$.
In each case, we have
$$\sigma(k)\ne k, \quad \sigma(\sigma(k))=k, \quad \mbox{and} \quad \sigma(2^ik)=2^i\sigma(k)$$
(the last equality is derived from
$\Zech(2k)
= \log_{\xi}(1+\xi^{2k})
=  \log_{\xi}(1+\xi^{k})^2
= 2 \log_{\xi}(1+\xi^{k})
= 2\Zech(k)$).
If $\sigma(k)=2^ik$
for some~$i\in \{1,\ldots,{n-1}\}$,
then
$k = \sigma(\sigma(k)) = \sigma(2^ik)=2^i\sigma(k)=2^{2i}k$,
$k (2^{2i} -1) =0$,
and $2^{2i} - 1 \equiv 0 \bmod 2^n-1$,
which is impossible for odd~$n$.

(iii)
Let $\tau(k)$ denote $-\Zech(k)$ or $\Zech(-k)$.
In each case, we have
$\tau(\tau(\tau(k)))=k$
and $\tau(2^ik)=2^i\tau(k)$.
If $\tau(k)=2^ik$
for some~$i\in \{0,\ldots,{n-1}\}$,
then
$k = \tau(\tau(\tau(k))) = 2^{3i}k$,
$k (2^{3i} -1) =0$,
and $2^{3i} - 1 \equiv 0 \bmod 2^n-1$,
which implies $n\equiv 0 \bmod 3$ or~$i=0$.
In the case $i=0$, we have $3k=0$
(in particular, $n$ is even)
by Proposition~\ref{prop1}.
\end{proof}

\begin{corollary}\label{c:cosets6}
 If $n=6k+1$, then for any
 $k \in\ZZ_{2^n-1} \setminus\{0\}$
 the set~$\Gamma(k)$
 contains six elements from
 six distinct $2$-cyclotomic classes of the same size.
\end{corollary}

We denote
$$
\Cy \Gamma(k)
=
\bigcup_{s\in \Gamma(k)}\Cy(s). $$
By~$\Cos(t)$, we denote the union of all $2$-cyclotomic classes of the size~$t$, and by~$N_t$, the number of such cosets.
Note that $N_t>0$ only if $t$ divides~$n$.
The proof of the following proposition is the same as the proof
of \cite[Proposition\,1]{KLOS:Doob}, and we omit it here.

\begin{proposition}\label{N_s}
If $n$ is the product of primes where each prime is of form $6k+1$, then $N_t$ is divisible by~$18$ for every $t>1$.
\end{proposition}

\begin{conjecture}\label{conj2}
If $n$ is the product of primes where each prime is of form $6k+1$, then for every integer~$t>1$ there are $\frac{N_t}{18}$ pairs $(a_i,b_i)$, $i=1,\ldots, \frac{N_t}{18}$,
such that
$$
\Cos(t) =
\dot\bigcup_i
\big(\Cy\Gamma(a_i)
\,\dot\cup\,
\Cy\Gamma(b_i)
\,\dot\cup\,
\Cy\Gamma(a_i+b_i)\big) .$$
\end{conjecture}

\begin{remark}
 Note that
there are nonprime~$n$ of form~$6k+1$, for example, $n=25$,
for which
the conclusion
of Conjecture~\ref{conj2}
does not hold,
because the number of $2$-cyclotomic classes is not divisible by~$18$.
\end{remark}

\begin{proposition}
 Assume that $n$ is the product of primes of form $6k+1$.
 If the conclusion of Conjecture~\ref{conj2}
 holds, then the conclusion
 of Conjecture~\ref{conj1} also holds and there is a
 balanced $\param{n}{2}$-design
 invariant under the action
 of~$\FF_{2^n}^\MLT$
 and the action
 of~$\Aut(\FF_{2^n})$.
\end{proposition}
\begin{proof}
Straightforwardly,
if the conclusion
 of Conjecture~\ref{conj2}
 holds with the collection of
pairs $( a_{i}, b_{i} )$,
 $i=1,\ldots, \frac{N_t}{18}$ of each divisor $t>1$ of $n$, then the conclusion
 of Conjecture~\ref{conj1}
 holds with
 the collection of corresponding triples
 $(2^ja_{i},2^jb_{i}, -2^ja_{i}{-}2^jb_{i})$,
 $j=0,\ldots,t-1$, $t>1$, $t|n$.
 Since the resulting $\param{n}{2}$-design
 is invariant under $\FF_{2^n}^\MLT$,
 which acts transitively on points,
 it is balanced by the definition.
\end{proof}

\begin{theorem}\label{th:n7n13n19}
 For $n=7,13,19$,
 there is a
 balanced $\param{n}{2}$-design
 invariant under the action
 of~$\FF_{2^n}^\MLT$
 and under the action
 of~$\Aut(\FF_{2^n})$.
\end{theorem}
\begin{proof}
By solving the corresponding exact covering problem,
 it has been found that
  Conjecture~\ref{conj2}
  holds with
  \begin{itemize}
   \item $\FF_{2^{7}} =
  \FF_2[x] /
  \langle
  x^7+x+1
  \rangle$
  and $(a_1,b_1) = (1,9)$.
   \item $\FF_{2^{13}} =
  \FF_2[x] /
  \langle
  x^{13}+x^4+x^3+ x+1
  \rangle$
  and $(a_i,b_i)$
  from the following list:
 $(3, 3543)$,
 $(5, 1826)$,
 $(11, 1205)$,
 $(15, 1956)$,
 $(21, 2065)$,
 $(23, 3609)$,
 $(27, 883)$,
 $(29, 258)$,
 $(35, 1574)$,
 $(37, 158)$,
 $(39, 4008)$,
 $(43, 1384)$,
 $(59, 1463)$,
 $(67, 1877)$,
 $(71, 198)$,
 $(77, 725)$,
 $(81, 2643)$,
 $(83, 1051)$,
 $(93, 2760)$,
 $(107, 3858)$,
 $(113, 3042)$,
 $(115, 2882)$,
 $(121, 3629)$,
 $(135, 1502)$,
 $(139, 963)$,
 $(143, 1479)$,
 $(149, 586)$,
 $(151, 2486)$,
 $(167, 1469)$,
 $(171, 483)$,
 $(181, 1489)$,
 $(211, 2265)$,
 $(223, 3211)$,
 $(241, 882)$,
 $(273, 2984)$.
  \item $\FF_{2^{19}} =
  \FF_2[x] /
  \langle
  x^{19}+x^5+x^2+x+1
  \rangle$
  and $(a_i,b_i)$
  from the following dataset:
  \url{https://doi.org/10.17632/rskwsgp5y4.1}.
  \end{itemize}
The corresponding triangle design is constructed as the union
of the orbits of the triangles
 $\tri {1}{\xi^{a_i}}{\xi^{-b_i}}$ under the action of the semi-direct product
of the groups~$\FF_{2^n}^\MLT$ and~$\Aut(\FF_{2^n})$.
\end{proof}

\section{Theoretical constructions}\label{s:th}
In this section, we propose two theoretical constructions
of triangle designs over~$\FF_2$. The first one is recursive;
it produces $\param{n}{2}$-designs
for any $n\equiv 0,1 \bmod 6$ and balanced $\param{n}{2}$-designs
for any $n\equiv 1 \bmod 6$.
The second construction is direct; it uses a~$\param{6}{2}$-design, as well as a~$\paramg{12}{6}{2}$-design,
to produce a balanced $\paramg{n}{6}{2}$-design first,
and then an unbalanced $\param{n}{2}$-design, $n\equiv 0 \bmod 6$.

\subsection{Recursive construction of  triangle designs}\label{s:rec}
The next theorem
does not give new existence results,
compared to the known results,
but we use it as an intermediate step
to construct new balanced triangle designs
in the following Theorem~\ref{th:balanced-rec}.
\begin{theorem}\label{th:re_con}
If there are triangle designs $\TT_m$ in $\PG(\FF_2^m)$ and $\TT_n$ in $\PG(\FF_2^n)$ and $mn$ is even,
then there is a triangle design in $\PG(\FF_2^{m+n})$.
\end{theorem}
\begin{proof}
Denote the sets of all lines of
$\PG(\FF_2^{m})$,
$\PG(\FF_2^{n})$
and $\PG(\FF_2^{m+n})$
by $\UU_m$, $\UU_n$ and $\UU_{m+n}$,
respectively.
We have,
$$|\UU_m|=\frac{(2^m-1)(2^m-2)}{6}=3|\TT_m|,$$
$$|\UU_n|=\frac{(2^n-1)(2^n-2)}{6}=3|\TT_n|, \quad\mbox{and}$$
$$|\UU_{m+n}|=\frac{(2^{m+n}-1)(2^{m+n}-2)}{6}.$$
We consider $\FF_2^{m+n}$ as the direct product $\FF_2^{m} \times \FF_2^{n}$. The lines in
$\UU_{m+n}$ are divided into the following types:
\begin{itemize}
\item[(a)] $\{(\vc{x},\vc{0}),
             (\vc{y},\vc{0}),
             (\vc{z},\vc{0})\}$,
    where $ \{ \vc x , \vc y , \vc z \} \in \UU_{m} $. The number of lines in this type is $|\UU_m|$;
\item[(b)] $\{(\vc{0},\vc{u}),
             (\vc{0},\vc{v}),
             (\vc{0},\vc{w})\}$,
    where $ \{ \vc u , \vc v , \vc w \}  \in \UU_{n} $. The number of lines in this type is $|\UU_n|$;
\item[(c)] $\{(\vc{x},\vc{u}),
             (\vc{y},\vc{v}),
             (\vc{z},\vc{w})\}$,
    where $ \{ \vc x , \vc y , \vc z \}  \in \UU_{m} $, $ \{ \vc u , \vc v , \vc w \}  \in \UU_{n} $. The number of lines in this type is $6|\UU_m||\UU_n|$;
\item[(d)] $\{ (\vc{x},\vc{0}),
          (\vc{y},\vc{v}),
          (\vc{z},\vc{v}) \}$,
    where $ \{ \vc x , \vc y , \vc z \}  \in \UU_{m} $
    and $\vc{0}\ne \vc{v} \in \FF_2^{n} $. The number of lines in this type is $3(2^n-1)|\UU_m|$;
\item[(e)] $\{ (\vc{0},\vc{u}),
          (\vc{y},\vc{v}),
          (\vc{y},\vc{w}) \}$,
    where $ \{ \vc u , \vc v , \vc w \}  \in \UU_{n} $
    and $ \vc{0}\ne\vc{y} \in \FF_2^{m} $. The number of lines in this type is $3(2^m-1)|\UU_n|$;
\item[(f)] $\{ (\vc{0},\vc{v}),
          (\vc{y},\vc{0}),
          (\vc{y},\vc{v}) \}$,
    where $ \vc{0}\ne\vc{y} \in \FF_2^{m} $
    and $\vc{0}\ne \vc{v} \in \FF_2^{n} $. The number of lines in this type is $(2^m-1)(2^n-1)$.
\end{itemize}

It is easy to check that the total number of lines in these six types is $|\UU_{m+n}|$. Assume w.l.o.g. that $n$ is even. In this case, there is a spread $S$ of the set of points
of~$\PG(\FF_2^n)$ into lines. We will construct a triangle design $\TT_{m+n} = \TT_m \times_S \TT_n$ based
 on the triangle designs~$\TT_m$
 and~$\TT_n$
 and spread~$S$, where the size of $S$ is $|S|=(2^n-1)/3$ and the size of $\TT_{m+n}$ is $|\TT_{m+n}|=|\UU_{m+n}|/3$.
 We include in~$\TT_{m+n}$ the following triangles:
\begin{itemize}
    \item[(A)] the triangle
    $$\big\{ a \times \{\vc{0}\}, b \times \{\vc{0}\}, c \times \{\vc{0}\} \big\}$$
    for each triangle $\{ a, b, c \}$ in~$\TT_m$. The number of triangles in this type is $|\UU_m|/3$;
    \item[(B)] the triangle
    $$\big\{\{\vc{0}\} \times a , \{\vc{0}\} \times b,  \{\vc{0}\}\times c \big\}$$
    for each triangle $\{ a, b, c \}$ in~$\TT_n$. The number of triangles in this type is $|\UU_n|/3$;
    \item[(C)] the triangle
    $$\big\{\{(\vc x,\vc u),(\vc y',\vc v),(\vc z,\vc w)\},
       \{(\vc z,\vc w),(\vc x',\vc u),(\vc y,\vc v)\},
       \{(\vc y,\vc v),(\vc z',\vc w),(\vc x,\vc u)\}
       \big\},$$
    for each triangle
    $\{\{\vc x,\vc y',\vc z\},
       \{\vc z,\vc x',\vc y\},
       \{\vc y,\vc z',\vc x\}
       \}$ in~$\TT_m$
    and each ordered triple $(\vc u,\vc v, \vc w)$ of vectors of~$\FF_2^n$ forming a line
    (each line corresponds to $3!=6$ triples). The number of triangles in this type is $2|\UU_m||\UU_n|$;
\item[(D)]
the triangle
    $$\big\{\{(\vc x,\vc v),(\vc y,\vc 0),(\vc z,\vc v)\},
       \{(\vc z,\vc v),(\vc x,\vc 0),(\vc y,\vc v)\},
       \{(\vc y,\vc v),(\vc z,\vc 0),(\vc x,\vc v)\}
       \big\},$$
for each line
$\{\vc x,\vc y, \vc z\}$ in~$\UU_m$
and each nonzero $\vc v$ in~$\FF_2^n$. The number of triangles in this type is $(2^n-1)|\UU_m|$;
\item[(E)]
the triangle
    $$\big\{\{(\vc y,\vc u),(\vc 0,\vc v),(\vc y,\vc w)\},
       \{(\vc y,\vc w),(\vc 0,\vc u),(\vc y,\vc v)\},
       \{(\vc y,\vc v),(\vc 0,\vc w),(\vc y,\vc u)\}
       \big\},$$
for each line
$\{\vc u,\vc v, \vc w\}$ in~$\UU_n\backslash S$ and each nonzero $\vc y$ in~$\FF_2^m$. The number of triangles in this type is $(2^m-1)(|\UU_n|-|S|)$;
\item[(F)]
the two triangles
\begin{align*}
   \big\{\{(\vc y,\vc u),(\vc 0,\vc u),(\vc y,\vc 0)\},
        \{(\vc y,\vc 0),(\vc y,\vc v),(\vc 0,\vc v)\},
        \{(\vc 0,\vc v),(\vc y,\vc w),(\vc y,\vc u)\}
        \big\},\\
    \big\{\{(\vc y,\vc w),(\vc y,\vc 0),(\vc 0,\vc w)\},
       \{(\vc 0,\vc w),(\vc y,\vc u),(\vc y,\vc v)\},
       \{(\vc y,\vc v),(\vc 0,\vc u),(\vc y,\vc w)\}
       \big\},
\end{align*}
for each line
$\{\vc u,\vc v, \vc w\}$ from~$S$ and each nonzero $\vc y$ in~$\FF_2^m$. The number of triangles in this type is $2(2^m-1)|S|$.
\end{itemize}
It is easy to check that the total number of triangles in these six types is $|\UU_{m+n}|/3=|\TT_{m+n}|$. It is now straightforward to check
that
\begin{itemize}
\item each line of type (a) belongs to a triangle of type (A); the number of lines in this case is $|\UU_m|$;
\item each line of type (b) belongs to a triangle of type (B); the number of lines in this case is $|\UU_n|$;
\item each line of type (c) belongs to a triangle of type (C); the number of lines in this case is $6|\UU_m||\UU_n|$;
\item each line of type (d) belongs to a triangle of type (D); the number of lines in this case is $3(2^n-1)|\UU_m|$;
\item each line of type (e) where $\{\vc u,\vc v, \vc w\}\not\in S$ belongs to a triangle of type (E); the number of lines in this case is $3(2^m-1)(|\UU_n|-|S|)$;
\item each line of type (e) where $\{\vc u,\vc v, \vc w\}\in S$ belongs to a triangle of type (F); the number of lines in this case is $3(2^m-1)|S|$;
\item each line of type (f) belongs to a triangle of type (F); the number of lines in this case is $(2^m-1)(2^n-1)$.
\end{itemize}
By counting, every line in $\UU_{m+n}$ belongs to exactly one triangle in $\TT_{m+n}$, which means that $\TT_{m+n}$ is a triangle design in $\FF_2^{m+n}$.
\end{proof}

Unfortunately, this recursive construction cannot start
from trivial values,
and we need at least a $\param{6}{2}$-design
to start the recursive process.
However,
the trivial $\param{1}{2}$-design,
which consists of the empty triangle set,
is useful  to reach odd values of dimension
from even ones, for example, $\param{7}{2}$
from~$\param{6}{2}$.

\begin{corollary}
 By the construction in Theorem~\ref{th:re_con},
 one can recursively construct an
 $\param{n}{2}$-design for any
 $n\equiv 1,0 \bmod 6$ from the
 $\param{6}{2}$-design constructed in Theorem~\ref{th:6_2}
 and the trivial $\param{1}{2}$-design.
\end{corollary}

We can observe that the construction
in Theorem~\ref{th:re_con} produces an unbalanced design even if
$m+n$ is odd,
because $m$ or $n$ is even and hence
$\TT_m$ or $\TT_n$ is unbalanced.
In the rest of this section,
we modify it to construct a balanced
triangle design.

\begin{theorem}\label{th:balanced-rec}
If there is a balanced triangle design $\TT_m$ in $\PG(\FF_2^m)$,
$m\ge 7$,
then there is a balanced triangle design in $\PG(\FF_2^{m+6})$.
\end{theorem}
\begin{proof}
We consider the construction in the proof
of the previous theorem with
$n=6$ and $\TT_6$ from Section~\ref{sec:6_2}
and modify the construction to make it balanced.
For a collection $\TT$ of triangles,
we will say that it gives the \emph{charge}
$a-b$ to a point~$\vc{x}$ if $\vc{x}$ belongs to exactly~$a$
triangles of~$\TT$ as a \emph{corner}
(i.e., lies in two lines of the triangle)
and to exactly~$b$
triangles of~$\TT$ as a \emph{non-corner}
(i.e., lies in only one line of the triangle).
We note that a (group divisible)
triangle design is balanced if and only if
it gives the zero charge to each point.
We are now going to consider each of the parts (A)--(F) of the construction in Theorem~\ref{th:re_con} separately and, if necessary, modify them to make the result balanced.

(A) This part is based on the balances design $\TT_m$,
gives the charge $0$ to all points
of the form $(\vc{x},\vc{0})$,
and does not touch any other point.
We keep this part unchanged.

(B) This part is based on the unbalanced design $\TT_6$;
it gives some nonzero charge to each point of the form $(\vc{0},\vc{v})$,
depending on~$\vc{v}$.
From the construction of~$\TT_6$
in Section~\ref{sec:6_2}, we see the following (where $\FF_2^6$ is associated with
$\FF_2 \times \FF_{2^5}$):
if $\vc{v}=(1,\vc{0})$,
then the charge is~$-31$;
if $\vc{v}=(0,\xi^i)$, $i\in\{0,\ldots,30\}$,
then the charge is~$+2$;
if $\vc{v}=(1,\xi^i)$, $i\in\{0,\ldots,30\}$,
then the charge is~$-1$.
We keep this part unchanged,
but nautralize the charge at step~(F).

(C) For each triangle~$T$ in~$\TT_m$
and each line~$L$ in~$\UU_n$, step (C) makes
$6$ triangles, where $9$ points get the charge~$+2$
(namely, the points of form~$(\vc{x},\vc{v}$),
where~$\vc{x}$ is a corner of~$T$ and $\vc{v}\in L$),
 and $9$ points get the charge~$-2$
(the points of form~$(\vc{x}',\vc{v}$),
where~$\vc{x}'$ is a non-corner of~$T$ and $\vc{v}\in L$). However, when $T$ runs
over all triangles of~$\TT_m$,
the charge of each point sums to~$0$ because $\TT_m$ is balanced.
So, we keep part~(C) unchanged.

(D) This part, as it is in Theorem~\ref{th:re_con},
is not good for us because
it gives only negative charges
to points of form~$(\vc{x},\vc{0})$.
Instead, we group the lines in~$\UU_m$ into triangles
$T=\{\{\vc x,\vc y',\vc z\},
       \{\vc z,\vc x',\vc y\},
       \{\vc y,\vc z',\vc x\}
       \}$ of~$\TT_m$ and construct the triangles
\begin{align*}\big\{\{(\vc x,\vc 0),(\vc y',\vc v),(\vc z,\vc v)\},
\{(\vc z,\vc v),(\vc x',\vc 0),(\vc y,\vc v)\},
       \{(\vc y,\vc v),(\vc z',\vc v),(\vc x,\vc 0)\}
       \big\},
       \\
\big\{\{(\vc x,\vc v),(\vc y',\vc 0),(\vc z,\vc v)\},
       \{(\vc z,\vc v),(\vc x',\vc v),(\vc y,\vc 0)\},
       \{(\vc y,\vc 0),(\vc z',\vc v),(\vc x,\vc v)\}
       \big\},
       \\
       \big\{\{(\vc x,\vc v),(\vc y',\vc v),(\vc z,\vc 0)\},
       \{(\vc z,\vc 0),(\vc x',\vc v),(\vc y,\vc v)\},
       \{(\vc y,\vc v),(\vc z',\vc 0),(\vc x,\vc v)\}
       \big\}.\end{align*}
Clearly, the lines involved are the same as those
in part~(D) in Theorem~\ref{th:re_con},
for the corresponding three lines of~$T$,
but the charges are different.
The points
$(\vc x ,\vc 0)$, $(\vc y ,\vc 0)$, $(\vc z ,\vc 0)$ get $+1$;
$(\vc x',\vc 0)$, $(\vc y',\vc 0)$, $(\vc z',\vc 0)$ get $-1$;
$(\vc x ,\vc v)$, $(\vc y ,\vc v)$, $(\vc z ,\vc v)$ get $+2$;
$(\vc x',\vc v)$, $(\vc y',\vc v)$, $(\vc z',\vc v)$ get $-2$.
When $T$ runs
over all triangles of~$\TT_m$,
the charge of each of these points sums to~$0$
because $\TT_m$ is balanced.

(E) This part is replaced similarly to (D), with the only difference
that now we group into triangles not all lines in~$\UU_6$,
but only the lines from~$\UU_6\backslash S$. For
this grouping, we use
the balanced $\paramg{6}{2}{2}$-design $\Gamma_{6,2}=(\FF_2^6,S,\BB)$
constructed in Theorem~\ref{th:gd6-2}. It should be noted
that the spread $S=S_{\vc{y}}$ now depends on~$\vc{y}$
(see the next part),
but all these spreads are equivalent;
so, $\Gamma_{6,2}$ also depends on~$\vc{y}$,
and all these group divisible triangle designs
can be constructed as equivalent.

(F) This part of construction in Theorem~\ref{th:re_con},
for each line $\{\vc u,\vc v, \vc w\}$ in the spread~$S$, and each~$\vc y$ in~$\FF_2^m$,
charges $(\vc 0,\vc u)$ by~$-2$
and charges each of $(\vc 0,\vc v)$, $(\vc 0,\vc w)$ by~$+1$. The vectors $(\vc y,\vc u)$, $(\vc y,\vc v)$, $(\vc y,\vc w)$, and $(\vc y,\vc 0)$ get the charge~$0$;
the other vectors are untouched.
Since $\vc u$, $\vc v$, $\vc w$ play the same role in the construction,
we are free to choose which of them gets~$-2$.
In particular, for three different~$\vc y$
from~$\FF_2^m$,
we can make this choice different
for the resulting charge to sum to~$0$.
We do this for $2^m-32$ values of~$\vc y$ (since $\TT_m$ is balanced, $m$ is odd and $2^m-32$ is divisible by~$3$).

The remaining $31$ values of~$\vc y$,
call them~$\vc y_t$, $i=0,1,\ldots,30$,
are intended to cancel the nonzero charge
from part~(B) as follows.
It is not difficult
to see that the spread~$S$ in~$\FF_2^6$
consists of $21$ lines,
one line of form
$\{(1,\vc{0}),(0,\xi^i),(1,\xi^i)\}$,
$15$ lines of form
$\{(0,\xi^i),(1,\xi^j),(1,\xi^k)\}$,
and $5$ lines of form
$\{(0,\xi^i),(0,\xi^j),(0,\xi^k)\}$.
For $\vc y=\vc y_0$ and each of these lines,
we construct two triangles as in part~(F) of the proof of Theorem~\ref{th:re_con},
caring about the charges.
\begin{itemize}
 \item[(i)]  For the line of form
$\{(1,\vc{0}),(0,\xi^i),(1,\xi^i)\}$,
we make $(\vc 0,(1,\vc{0}))$
and $(\vc 0,(1,\xi^i))$
to receive the charge~$+1$, $(\vc 0,(1,\xi^i))$ to receive the charge~$-2$.
 \item [(ii)]
For the~$15$ lines of form
$\{(0,\xi^i),(1,\xi^j),(1,\xi^k)\}$,
we make $(\vc 0,(0,\xi^i))$ to receive the charge~$-2$ in $5$ cases,  and charge~$-2$ in the remaining $10$ cases. The vectors
$(\vc 0,(1,\xi^i))$, $(\vc 0,(1,\xi^k))$
get the charges $+1$, $+1$ in the first $5$ cases
and $+1$, $-2$ in the last $10$ cases.
\item[(iii)] For each of the remaining $5$ lines, the choice is not important, because the first coordinate is $0$ for all points of the line.
\end{itemize}
For $\vc y =\vc y_t$,
we choose the same charges as for $\vc y_t$,
but multiply the last $\FF_{2^5}$-component
of each vector by~$\xi^t$ (this changes the spread~$S$ as well, which is not a problem).
With such ``rotation'', we see that the
summary charge of all vectors from the part of construction described in~(ii) and~(iii) is~$0$.
The charge from~(i) sums to $+31$ for $(\vc 0,(1,\vc 0))$, is~$+1$ for each
$(\vc 0,(1,\xi^i))$, and~$-2$
for each $(\vc 0,(1,\xi^i))$, $i=0,1,\ldots,30$.

Finally, we see that after the modification,
the construction covers the same lines as
the construction in Theorem~\ref{th:re_con};
each of parts (A), (C), (D), (E) does not
give any nonzero charge to any vector;
the charges from parts~(B) and~(F)
neutralize each other. Hence,
we have constructed a balanced triangle design.
\end{proof}

\subsection
[Balanced gd(6k,6){2}-designs and related
unbalanced (6k){2}-designs]
{Balanced \boldmath$\paramg{6k}{6}{2}$-designs
and \boldmath$\paramg{6k}{2}{2}$-designs}
\label{sec:6k}

From the balanced $\paramg{12}{6}{2}$-design found in Section~\ref{sec:gd12_6_2},
we can construct an infinite series of
balanced group-divisible triangle designs.

\begin{theorem}\label{th:6k,6}
 For every positive integer~$k$,
 there exists a balanced $\paramg{6k}{6}{2}$-design
 $(\VV,\GG,\BB)$, where
 $\VV$ is a $k$-dimensional space over~$\FF_{2^6}$
 and~$\GG$ is the set of its $1$-dimensional
 $\FF_{2^6}$-subspaces.
\end{theorem}
\begin{proof}
 Let $k\geq 2$ be a positive integer.
Let $\VV=\FF_{2^{6k}}$.
Then $\VV$ can be treated as
a $k$-dimensional vector space over~$\FF_{2^6}$;
also, when we discuss lines and triangles in $\PG(\VV)$ it is treated
as a $6k$-dimensional vector space over~$\FF_{2}$.
Denote by $\GG$ and $\DD$ the sets of all $1$- and  $2$-dimensional subspaces of~$\VV$ over~$\FF_{2^6}$, respectively. By Theorem~\ref{th:gd12-6},
  for each~$D$ from~$\DD$,
  there is a set~$\BB_D$ of triangles in~$D$
  such that $(D,\GG_D,\BB_D)$ is a
  balanced $\paramg{12}{6}{2}$-design,
  where~$\GG_D$, $\GG_D\subset\GG$, is the set of all $1$-dimensional $\FF_{2^6}$-subspaces of~$D$. By the definition,
  $
  (\VV,\GG,\bigcup_{D\in \DD}\BB_D)
  $ is a $\paramg{6k}{6}{2}$-design.
  Indeed, if a line~$L$
  lies in a $1$-dimensional
  $\FF_{2^6}$-subspace,
  then it does not belong to any triangle from~$\BB_D$
  for any~$D$ from~$\DD$.
  Otherwise, $L$ spans a unique
  $2$-dimensional $\FF_{2^6}$-subspace~$D$
  and belongs to a unique triangle from~$\BB_D$.

  It remains to note that the group-divisible
  triangle design $
  (\VV,\GG,\bigcup_{D \in \DD}\BB_D)
  $ is balanced because each
  of the designs $(D,\GG_D,\BB_D)$ is balanced and each point
  belongs to a constant number of $2$-dimensional
  $\FF_{2^6}$-subspaces~$D$.
\end{proof}

It is easy to see that the proof of
Theorem~\ref{th:6k,6} generalizes to the following.
\begin{lemma}\label{l}
Let $\VV$ be a $2$-dimensional vector space over~$\FF_{q^m}$. Let~$\GG$ be the set of all $1$-dimensional subspaces of~$\VV$ over~$\FF_{q^m}$.
Assume that there exists a (balanced)
$\paramg{2m}{m}{q}$-design
$(\VV,\GG,\BB)$.
Then there exists a (balanced)
$\paramg{km}{m}{q}$-design for all~$k$.
\end{lemma}

The following lemma is straightforward.
\begin{lemma}\label{l:gd-t}
 If there is
 a $\paramg{n}{m}{q}$-design
 and
 a $\paramg{m}{t}{q}$-design, then there is a
 $\paramg{n}{t}{q}$-design (in particular, a $\param{n}{q}$-design if $t=1$).
 If there is
 a balanced $\paramg{n}{m}{q}$-design
 and
 a balanced $\paramg{m}{t}{q}$-design, then there is a
 balanced $\paramg{n}{t}{q}$-design.
\end{lemma}
By Theorems~\ref{th:gd6-2} and~\ref{th:6k,6} and Lemma~\ref{l:gd-t}, we get the following.
\begin{corollary}\label{cor:gd6k-2}
 For every positive integer $k$,
 there exists a balanced $\paramg{6k}{2}{2}$-design.
\end{corollary}
By Theorems~\ref{th:6_2} and~\ref{th:6k,6} and Lemma~\ref{l:gd-t}, we get the following (which, actually, gives only a new construction for the known existence result).
\begin{corollary}\label{cor:6k}
 For every positive integer $k$,
 there exists a $\param{6k}{2}$-design.
\end{corollary}

\section{Conclusion}\label{sec:7}
We have constructed balanced triangle designs over~$\FF_{2}$
for all~$n$ congruent to~$1$ modulo~$6$
and found three balanced triangle designs
that are symmetric with respect to the action of~$\FF_{2^n}^\MLT$
and~$\Aut(\FF_{2^n})$. We have conjectured that such
symmetric designs exist at least for each prime~$n$
congruent with~$0$ or~$1$ modulo~$6$
and related their existence with the existence
of special partitions of the set of $2$-cyclotomic
classes of~$\ZZ_{2^n}\setminus \{0\}$
into subsets of size~$18$ (Conjecture~\ref{conj2}).
Special partitions of the set of
$2$-cyclotomic classes into subsets
of size~$6$ satisfying some special conditions
were constructed in~\cite{KLOS:Doob}, related to a problem in coding theory.
 In that case, an algebraic construction was found, and we believe
 that Conjecture~\ref{conj2} has an analytic solution as well.

 Additionally, we have constructed
 an infinite class of
 balanced group divisible triangle designs, with groups of dimension~$6$ and~$2$.
 Taking into account the well-known relations between the set of lines of~$\PG(n-1,2)$ and combinatorial designs,
 the group divisible triangle designs constructed in Theorem~\ref{th:6k,6}
 and Corollary~\ref{cor:gd6k-2}
 can be treated as balanced triangle decompositions
 of group divisible designs GDD$(2^n-1,3,63)$ and GDD$(2^n-1,3,3)$, respectively,
 or generalizes Steiner systems
 GS$((2^n-1)/63,3,2;64)$
 and
 GS$((2^n-1)/3,2;4)$,
 in the sense of~\cite{Etzion:97}.

\subsection*{Acknowledgements}
The authors thank the reviewer for the suggestions, which helped to improve the presentation of the material,
and Professor Marco Buratti for his brilliant lecture at the Combinatorial Constructions Conference, Dubrovnik, April 7--13, 2024, which attracted the author's attention to Heffter's difference problem and related questions.

\subsection*{Data availability}
Data are available in Mendeley Data, \url{https://doi.org/10.17632/rskwsgp5y4.1}


\providecommand\href[2]{#2} \providecommand\url[1]{\href{#1}{#1}} \def\DOI#1{{\href{https://doi.org/#1}{https://doi.org/#1}}}\def\DOIURL#1#2{{\href{https://doi.org/#2}{https://doi.org/#1}}}

\end{document}